\documentclass[twosided,reqno]{amsart}

\usepackage{amsmath}
\usepackage{amssymb}
\usepackage{amsthm}

\theoremstyle{theorem}
\newtheorem{theorem}{\scshape Theorem }[section]

\theoremstyle{definition}

\newtheorem{remark}{\scshape Remark}

\numberwithin{equation}{section}

\begin{document}

\title{Higher-order Frobenius-Euler and poly-Bernoulli mixed type polynomials}

\author{Dae San Kim}
\address{Department of Mathematics, Sogang University, Seoul 121-742, Republic of Korea.}
\email{dskim@sogang.ac.kr}

\author{Taekyun Kim}
\address{Department of Mathematics, Kwangwon University, Seoul 139-701, Republic of Korea}
\email{tkkim@kw.ac.kr}

\maketitle

\begin{abstract}
In this paper, we consider higher-order Frobenius-Euler polynomials associated with poly-Bernoulli polynomials which are derived from polylogarithmic function.  These polynomials are
called higher-order Frobenius-Euler and poly-Bernoulli mixed type polynomials.  The purpose of this paper is to give various identities of those polynomials arising from umbral calculus.
\end{abstract}

\section{Introduction}

For $\lambda \in \mathbb{C}$ with $\lambda \neq 1$, the Frobenius-Euler polynomials of order $\alpha(\alpha \in \mathbb{R})$ are defined by the generating function to be

\begin{equation}\label{1}
\left(\frac{1-\lambda}{e^t-\lambda}\right)^{\alpha}e^{xt}=\sum_{n=0} ^{\infty}H_n ^{(\alpha)} (x|\lambda)\frac{t^n}{n!},{\text{ (see [1,6,7,13,14])}}.
\end{equation}

When $x=0$, $H_n ^{(\alpha)}(\lambda)=H_n  ^{(\alpha)}(0\vert \lambda)$ are called the Frobenius-Euler numbers of order $\alpha$. As is well known, the Bernoulli polynomials of order $\alpha$ are defined by the generating function to be

\begin{equation}\label{2}
\left(\frac{t}{e^t-1}\right)^{\alpha} e^{xt}=\sum_{n=0} ^{\infty}\mathbb{B}_n ^{(\alpha)} (x)\frac{t^n}{n!},{\text{ (see [4,5,9])}}.
\end{equation}

When $x=0$, $\mathbb{B}_n ^{(\alpha)}=\mathbb{B}_n  ^{(\alpha)}(x)$ is called the $n$-th Bernoulli number of order $\alpha$.
In the special case, $\alpha=1$, $\mathbb{B}_n ^{(1)}(x)=B_n(x)$ is called the $n$-th Bernoulli polynomial.  When $x=0, B_n=B_n(0)$ is
called the $n$-th ordinary Bernoulli number.  Finally, we recall that the Euler polynomials of order $\alpha$  are given by

\begin{equation}\label{3}
\left(\frac{2}{e^t+1}\right)^{\alpha} e^{xt}=\sum_{n=0} ^{\infty}E_n ^{(\alpha)} (x)\frac{t^n}{n!},{\text{ (see [2,3,8,10,15])}}.
\end{equation}

When $x=0$, $E_n ^{(\alpha)}=E_n  ^{(\alpha)}(0)$ is called the $n$-th Euler number of order $\alpha$.  In the special case, $\alpha =1$, $E_{n}^{(1)}(x)=E_{n}(x)$ is called
the  $n$-th ordinary Euler polynomial.  The classical polylogarithmic function $Li_{k}(x)$ is defined by

\begin{equation}\label{4}
Li_{k}(x)=\sum_{n=1}^{\infty}\frac{x^n}{n^k}, ~~~(k\in \mathbb{Z}),~~(see[5]).
\end{equation}

As is known, poly-Bernoulli polynomials are defined by the generating function to be

\begin{equation}\label{5}
\frac{Li_{k}(1-e^{-t})}{1-e^{-t}} e^{xt}=\sum_{n=0} ^{\infty}B_n ^{(k)} (x)\frac{t^n}{n!},{\text{ (cf. [5])}}.
\end{equation}

Let $\mathbb{C}$ be the complex number field and let ${\mathcal{F}}$ be the set of all formal power series in the variable $t$ over ${\mathbb{C}}$ with

\begin{equation}\label{6}
{\mathcal{F}}=\left\{ \left.f(t)=\sum_{k=0} ^{\infty} \frac{a_k}{k!} t^k~\right|~ a_k \in {\mathbb{C}} \right\}.
\end{equation}

Now, we use the notation $\mathbb{P}=\mathbb{C}[x]$.  In this paper, $\mathbb{P}^{*}$ will be denoted by the vector space of all linear functionals on ${\mathbb{P}}$.
Let us assume that $\left< L~ \vert ~p(x)\right>$ be the action of the linear functional $L$ on the polynomial $p(x)$, and we remind that the vector space operations on $\mathbb{P}^{*}$
are defined by $\left< L+M~ \vert ~p(x)\right>=\left< L~ \vert ~p(x)\right>+\left< M~ \vert ~p(x)\right>$, $\left< cL~ \vert ~p(x)\right>=c \left< L~ \vert ~p(x)\right>$, where $c$ is a complex constant in $\mathbb{C}$.  The formal power series

\begin{equation}\label{7}
f(t)=\sum_{k=0} ^{\infty} \frac{a_k}{k!} t^k~\in \mathcal{F},
\end{equation}

defines a linear functional on $\mathbb{P}$ by setting

\begin{equation}\label{8}
\left<f(t)|x^n \right>=a_n,~for~all~n\geq 0, {\text{(see [11,12])}}.
\end{equation}

From (\ref{7}) and (\ref{8}), we note that

\begin{equation}\label{9}
\left<t^k | x^n \right>=n! \delta_{n,k},~(see~ [11,12]),
\end{equation}

where $\delta_{n,k}$ is the Kronecker symbol.\\

Let us consider $f_{L}(t)=\sum_{k=0}^{\infty}\frac{\left<L | x^n \right>}{k!}t^k$.  Then we see that $\left<f_{L}(t)|x^n \right>=\left<L | x^n \right>$ and
so $L=f_{L}(t)$ as linear functionals.  The map $L \mapsto f_{L}(t)$ is a vector space isomorphism from $\mathbb{P}^{*}$ onto $\mathcal{F}$.  Henceforth, $\mathcal{F}$ will denote
both the algebra of formal power series in $t$ and the vector space of all linear functionals on $\mathbb{P}$, and so an element $f(t)$ of $\mathcal{F}$ will be thought of as both a formal power series and a linear functional (see[11]).  We shall call $\mathcal{F}$ the umbral algebra.  The umbral calculus is the study of umbral algebra.  The order $o(f(t))$ of a nonzero power series
$f(t)$ is the smallest integer $k$ for which the coefficient of $t^k$ does not vanish.  A series $f(t)$ is called a delta series if  $o(f(t))=1$, and an invertible seires if $o(f(t))=0$.  Let $f(t),g(t) \in \mathcal{F}$.  Then we have

\begin{equation}\label{10}
\left<f(t)g(t) | p(x) \right>=\left<f(t) | g(t)p(x) \right>= \left<g(t) | f(t)p(x) \right>,~(see~ [11]).
\end{equation}

For $f(t),g(t) \in \mathcal{F}$ with $o(f(t))=1,~o(g(t))=0$, there exists a unigue sequence $S_n(x)(deg~S_n(x)=n)$ such that
$\left<g(t)f(t)^k | S_n(x) \right>=n! \delta_{n,k}~for~n,k\geq0$. The sequence $S_n(x)$ is called the Sheffer sequence for $(g(t),f(t))$ which is
denoted by $S_n(x)\sim (g(t),f(t)),~(see~[11,12])$.  Let $f(t) \in \mathcal{F}$ and $p(t) \in \mathbb{P}$.  Then we have

\begin{equation}\label{11}
f(t)=\sum_{k=0} ^{\infty} \left<f(t)|x^k\right>\frac{t^k}{k!},~p(x)=\sum_{k=0} ^{\infty}\left<t^k|p(x)\right> \frac{x^k}{k!}.
\end{equation}

From (\ref{11}), we note that

\begin{equation}\label{12}
p^{(k)}(0)= \left<t^k|p(x)\right> =\left<1|p^{(k)}(x)\right>.
\end{equation}

By (\ref{12}), we get

\begin{equation}\label{13}
t^{k}p(x)= p^{(k)}(x)= \frac{d^kp(x)}{dx^k},~(see~[11,12]).
\end{equation}

From (\ref{13}), we easily derive the following equation:

\begin{equation}\label{14}
e^{yt}p(x)=p(x+y),~~\langle e^{yt}\vert p(x)\rangle =p(y).
\end{equation}

For $p(x) \in \mathbb{P},~f(t) \in \mathcal{F}$, it is known that

\begin{equation}\label{15}
\langle f(t)\vert xp(x)\rangle = \langle \partial_t f(t)\vert p(x)\rangle =\langle f'(t)\vert p(x)\rangle,~(see[11]).
\end{equation}

Let $S_n(x) \sim (g(t), f(t))$.  Then we have

\begin{equation}\label{16}
\frac{1}{g(\bar{f}(x))}e^{y\bar{f}(t)}= \sum_{n=0}^{\infty}S_n(y)\frac{t^n}{n!},~for~all~y \in \mathbb{C},
\end{equation}

where $\bar{f}(t)$ is the compositional inverse of $f(t)$ with $\bar{f}(f(t))=t$, and

\begin{equation}\label{17}
f(t)S_n(x)=nS_{n-1}(x),~~(see~[11,12]).
\end{equation}

The Stirling number of the second kind is defined by the generating function to be

\begin{equation}\label{18}
(e^t-1)^m=m!\sum_{l=m} ^{\infty}S_2(l,m)\frac{t^m}{m!},~(m \in \mathbb{Z}_{\geq 0}).
\end{equation}

For $S_n(x) \sim (g(t),t)$, it is well known that

\begin{equation}\label{19}
S_{n+1}(x)=(x-\frac{g'(t)}{g(t)})S_n(x),~~(n \geq 0),~~(see~[11,12]).
\end{equation}

Let $S_n(x) \sim (g(t),f(t))$, $r_n(x) \sim (h(t),l(t))$.  Then we have

\begin{equation}\label{20}
S_{n}(x)=\sum_{m=0}^{n}C_{n,m}r_m(x),
\end{equation}

where

\begin{equation}\label{21}
C_{n,m}=\frac{1}{m!} \langle \frac{h(\bar{f}(t))}{g(\bar{f}(t))}l(\bar{f}(t))^m \vert x^n \rangle, ~~(see~[11,12]).
\end{equation}

In this paper, we study higher-order Frobeniuns-Euler polynomials associated with poly-Bernoulli polynomials which are called higher-order Frobenius-Euler and poly-Beroulli mixed type polynomials.  The purpose of this paper is to give various identities of those polynomials arising from umbral calculus.

\section{Higher-order Frobenius-Euler polynomials associated poly-Bernoulli polynomials}

Let us consider the polynomials $T_n^{(r,k)}(x|\lambda)$, called higher-order Frobenius-Euler and poly-Bernoulli mixed type polynomials, as follows:

\begin{equation}\label{22}
\left(\frac{1- \lambda}{e^t-\lambda}\right)^r\frac{Li_{k}(1-e^{-t})}{1-e^{-t}}e^{xt}=\sum_{n=0}^{\infty}T_n^{(r,k)}(x|\lambda)\frac{t^n}{n!},
\end{equation}

where $\lambda \in \mathbb{C}$ with $\lambda \neq 1$, $r,k \in \mathbb{Z}$.\\

When $x=0$, $T_{n}^{(r,k)}(\lambda)=T_{n}^{(r,k)}(0|\lambda)$ is called the $n$-th higher-order Frobenius-Euler and poly-Bernoulli mixed type number.\\

From (\ref{16}) and (\ref{22}), we note that

\begin{equation}\label{23}
T_{n}^{(r,k)}(x|\lambda) \sim \left(g_{r,k}(t)=\left(\frac{e^t-\lambda}{1- \lambda}\right)^r \frac{1-e^{-t}}{Li_{k}(1-e^{-t})},t\right).
\end{equation}

By (\ref{17}) and (\ref{23}), we get

\begin{equation}\label{24}
t T_{n}^{(r,k)}(x|\lambda) =nT_{n-1}^{(r,k)}(x|\lambda).
\end{equation}

From (\ref{22}), we can easily derive the following equation :

\begin{equation}\label{25}
\begin{split}
T_{n}^{(r,k)}(x|\lambda) &=\sum_{l=0}^{n}\binom{n}{l}H_{n-l}^{(r)}(\lambda)B_{l}^{(k)}(x)\\
&=\sum_{l=0}^{n}\binom{n}{l}H_{n-l}^{(r)}(x| \lambda)B_{l}^{(k)}.\\
\end{split}
\end{equation}

By (\ref{16}) and (\ref{23}), we get

\begin{equation}\label{26}
T_{n}^{(r,k)}(x|\lambda)=\frac{1}{g_{r,k}(t)}x^n=\left(\frac{1- \lambda}{e^t-\lambda}\right)^r\frac{Li_{k}(1-e^{-t})}{1-e^{-t}}x^n.
\end{equation}

In [5], it is known that

\begin{equation}\label{27}
\frac{Li_{k}(1-e^{-t})}{1-e^{-t}}x^n=\sum_{m=0}^{n}\frac{1}{(m+1)^k}\sum_{j=0}^{m}(-1)^j \binom{m}{j}(x-j)^n.
\end{equation}

Thus, by (\ref{26}) and (\ref{27}), we get

\begin{equation}\label{28}
\begin{split}
T_{n}^{(r,k)}(x|\lambda) &=\left(\frac{1- \lambda}{e^t-\lambda}\right)^r \frac{Li_{k}(1-e^{-t})}{1-e^{-t}} x^n\\
&=\sum_{m=0}^{\infty}\frac{1}{(m+1)^k} \sum_{j=0}^{m}(-1)^j \binom{m}{j}\left(\frac{1- \lambda}{e^t-\lambda}\right)^r(x-j)^n\\
&=\sum_{m=0}^{n}\frac{1}{(m+1)^k}\sum_{j=0}^{m}(-1)^j \binom{m}{j}H_{n}^{(r)}(x-j|\lambda).\\
\end{split}
\end{equation}

By (\ref{1}), we easily see that

\begin{equation}\label{29}
H_{n}^{(r)}(x|\lambda)=\sum_{l=0}^{n}\binom{n}{l}H_{n-l}^{(r)}(\lambda)x^{l}.
\end{equation}

Therefore, by (\ref{28}) and (\ref{29}), we obtain the following theorem.

\begin{theorem}\label{thm1}
For $r,k \in \mathbb{Z},~n \geq 0$, we have
\begin{equation*}
\begin{split}
T_{n}^{(r,k)}(x|\lambda)&=\sum_{m=0}^{n}\frac{1}{(m+1)^k}\sum_{j=0}^{m}(-1)^j \binom{m}{j}\sum_{l=0}^{n}\binom{n}{l}H_{n-l}^{(r)}(\lambda)(x-j)^{l}\\
&=\sum_{l=0}^{n} \left\{\binom{n}{l}H_{n-l}^{(r)}(\lambda)\sum_{m=0}^{\infty}\frac{1}{(m+1)^k}\sum_{j=0}^{m}(-1)^j \binom{m}{j}\right\}(x-j)^{l}.\\
\end{split}
\end{equation*}
\end{theorem}

In [5], it is known that

\begin{equation}\label{30}
\frac{Li_{k}(1-e^{-t})}{1-e^{-t}}x^n=\sum_{j=0}^{n}\left\{\sum_{m=0}^{n-j}\frac{(-1)^{n-m-j}}{(m+1)^{k}}\binom{n}{j}m!S_{2}(n-j, m)\right\}x^{j}.
\end{equation}

By (\ref{26}) and (\ref{30}), we get

\begin{equation}\label{31}
\begin{split}
T_{n}^{(r,k)}(x|\lambda) &=\left(\frac{1- \lambda}{e^t-\lambda}\right)^r \frac{Li_{k}(1-e^{-t})}{1-e^{-t}} x^n\\
&=\sum_{j=0}^{n}\left\{\sum_{m=0}^{n-j}\frac{(-1)^{n-m-j}}{(m+1)^{k}}\binom{n}{j}m!S_{2}(n-j, m)\right\}\left(\frac{1- \lambda}{e^t-\lambda}\right)^r x^{j}\\
&=\sum_{j=0}^{n}\left\{\sum_{m=0}^{n-j}\frac{(-1)^{n-m-j}}{(m+1)^{k}}\binom{n}{j}m!S_{2}(n-j, m)\right\}H_{j}^{(r)}(x|\lambda).\\
\end{split}
\end{equation}

Therefore, by (\ref{29}) and (\ref{31}), we obtain the following theorem.

\begin{theorem}\label{thm2}
For $r,k \in \mathbb{Z},~n \in \mathbb{Z}_{\geq 0}$, we have
\begin{equation*}
T_{n}^{(r,k)}(x|\lambda)=\sum_{l=0}^{n}\left\{\sum_{j=l}^{n}\sum_{m=0}^{n-j}(-1)^{n-m-j}\binom{n}{j}\binom{j}{l}\frac{m!}{(m+1)^k}H_{j-l}^{(r)}(\lambda)S_{2}(n-j, m)\right\} x^{l}.
\end{equation*}
\end{theorem}

From (\ref{19}) and (\ref{23}), we have

\begin{equation}\label{32}
T_{n+1}^{(r,k)}(x|\lambda)=\left(x-\frac{g'_{r,k}(t)}{g_{r,k}(t)}\right)T_{n}^{(r,k)}(x|\lambda).
\end{equation}

Now, we note that,

\begin{equation}\label{33}
\begin{split}
\frac{g'_{r,k}(t)}{g_{r,k}(t)}&=\left( \log g_{r,k}(t) \right)'\\
&= \left( r \log (e^t- \lambda) -r \log(1-\lambda)+ \log(1-e^{-t})- \log Li_{k}(1-e^t) \right)'\\
&=r+ \frac{r\lambda }{e^t \lambda }+\left( \frac{t}{e^t-1} \right)\frac{Li_{k}(1-e^{-t})-Li_{k-1}(1-e^{-t})}{tLi_{k}(1-e^{-t})}.\\
\end{split}
\end{equation}

By (\ref{32}) and (\ref{33}), we get

\begin{equation}\label{34}
\begin{split}
T_{n+1}^{(r,k)}(x|\lambda)& = xT_{n}^{(r,k)}(x|\lambda)-rT_{n}^{(r,k)}(x|\lambda)-\frac{r\lambda}{1-\lambda}\left( \frac{1-\lambda}{e^{t}-\lambda} \right)^{r+1}\frac{Li_{k}(1-e^{-t})}{1-e^{-t}} x^{n}\\
& - \left( \frac{1-\lambda}{e^{t}-\lambda} \right)^{r}\frac{Li_{k}(1-e^{-t})-Li_{k-1}(1-e^{-t})}{t(1-e^{-t})}\left(\frac{t}{e^t-1}\right) x^n\\
&=(x-r)T_{n}^{(r,k)}(x|\lambda)-\frac{r\lambda}{1-\lambda}T_{n}^{(r+1,k)}(x|\lambda)\\
&-\sum_{l=0}^{n}\binom{n}{l}B_{n-l} \left( \frac{1-\lambda}{e^{t}-\lambda} \right)^{r}\frac{Li_{k}(1-e^{-t})-Li_{k-1}(1-e^{-t})}{t(1-e^{-t})} x^l.\\
\end{split}
\end{equation}

It is easy to show that

\begin{equation}\label{35}
\begin{split}
\frac{Li_{k}(1-e^{-t})-Li_{k-1}(1-e^{-t})}{1-e^{-t}}&= \frac{1}{1-e^{-t}}\sum_{n=1}^{\infty}\{\frac{(1-e^{-t})^{n}}{n^{k}}-\frac{(1-e^{-t})^{n}}{n^{k-1}} \}\\
&=\left( \frac{1-e^{-t}}{2^{k}} - \frac{1-e^{-t}}{2^{k-1}}\right)+ \cdots \\
&=\left( \frac{1}{2^{k}} - \frac{1}{2^{k-1}}\right)t+ \cdots . \\
\end{split}
\end{equation}

For any delta series $f(t)$, we have

\begin{equation}\label{36}
\frac{f(t)}{t}x^{n}=f(t)\frac{1}{n+1}x^{n+1}.
\end{equation}

Thus, by (\ref{34}), (\ref{35}) and (\ref{36}), we get

\begin{equation}\label{37}
\begin{split}
T_{n+1}^{(r,k)}(x|\lambda)&=(x-r)T_{n}^{(r,k)}(x|\lambda)-\frac{r\lambda}{1-\lambda}T_{n}^{(r+1,k)}(x|\lambda)\\
&-\sum_{l=0}^{n}\binom{n}{l}B_{n-l}\frac{1}{l+1} \left( \frac{1-\lambda}{e^{t}-\lambda} \right)^{r} \frac{Li_{k}(1-e^{-t})-Li_{k-1}(1-e^{-t})}{1-e^{-t}} x^{l+1}\\
&=(x-r)T_{n}^{(r,k)}(x|\lambda)-\frac{r\lambda}{1-\lambda}T_{n}^{(r+1,k)}(x|\lambda)\\
&\qquad \qquad\qquad -\sum_{l=0}^{n}\frac{\binom{n}{l}}{l+1}B_{n-l}\{ T_{l+1}^{(r,k)}(x|\lambda)-T_{l+1}^{(r,k-1)}(x|\lambda)\}\\
&=(x-r)T_{n}^{(r,k)}(x|\lambda)-\frac{r\lambda}{1-\lambda}T_{n}^{(r+1,k)}(x|\lambda)\\
&\qquad \qquad\qquad -\frac{1}{n+1}\sum_{l=1}^{n+1}\binom{n+1}{l}B_{n+1-l}\{ T_{l}^{(r,k)}(x|\lambda)-T_{l}^{(r,k-1)}(x|\lambda)\}\\
&=(x-r)T_{n}^{(r,k)}(x|\lambda)-\frac{r\lambda}{1-\lambda}T_{n}^{(r+1,k)}(x|\lambda)\\
&\qquad \qquad\qquad -\frac{1}{n+1}\sum_{l=0}^{n+1}\binom{n+1}{l}B_{n+1-l}\{ T_{l}^{(r,k)}(x|\lambda)-T_{l}^{(r,k-1)}(x|\lambda)\}\\
&=(x-r)T_{n}^{(r,k)}(x|\lambda)-\frac{r\lambda}{1-\lambda}T_{n}^{(r+1,k)}(x|\lambda)\\
&\qquad \qquad\qquad -\frac{1}{n+1}\sum_{l=0}^{n+1}\binom{n+1}{l}B_{l}\{ T_{n+1-l}^{(r,k)}(x|\lambda)-T_{n+1-l}^{(r,k-1)}(x|\lambda)\}.\\
\end{split}
\end{equation}

Therefore, by (\ref{37}), we obtain the following theorem.

\begin{theorem}\label{thm3}
For $r,k \in \mathbb{Z}, ~n \in \mathbb{Z}_{\geq 0}$, we have
\begin{equation*}
\begin{split}
T_{n+1}^{(r,k)}(x|\lambda)
&=(x-r)T_{n}^{(r,k)}(x|\lambda)-\frac{r\lambda}{1-\lambda}T_{n}^{(r+1,k)}(x|\lambda)\\
&\qquad \qquad\qquad -\frac{1}{n+1}\sum_{l=0}^{n+1}\binom{n+1}{l}B_{l}\{ T_{n+1-l}^{(r,k)}(x|\lambda)-T_{n+1-l}^{(r,k-1)}(x|\lambda)\}.\\
\end{split}
\end{equation*}
\end{theorem}

\begin{remark}
If $r=0$, then we have
\begin{equation}\label{38}
\sum_{n=0}^{\infty}B_{n}^{(k)}(x)\frac{t^{n}}{n!}=\frac{Li_{k}(1-e^{-t})}{(1-e^{-t})}e^{xt}=\sum_{n=0}^{\infty} T_{n}^{(0,k)}(x|\lambda) \frac{t^n}{n!}.
\end{equation}
\end{remark}
Thus, by (\ref{38}), we get $B_{n}^{(k)}(x)=T_{n}^{(0,k)}(x|\lambda)$.

From (\ref{25}), we have

\begin{equation}\label{39}
\begin{split}
txT_{n}^{(r,k)}(x|\lambda)&=t\left( x\sum_{l=0}^{n}\binom{n}{l}H_{n-l}^{(r)}(\lambda) B_{l}^{(k)}(x) \right)\\
&=\sum_{l=0}^{n}\binom{n}{l}H_{n-l}^{(r)}(\lambda)\{lxB_{l-1}^{(k)}(x)+B_{l}^{(k)}(x)\}\\
&=nx\sum_{l=0}^{n-1}\binom{n-1}{l}H_{n-1-l}^{(r)}(\lambda)B_{l}^{(k)}(x)+\sum_{l=0}^{n}\binom{n}{l}H_{n-l}^{(r)}(\lambda)B_{l}^{(k)}(x)\\
&=nxT_{n-1}^{(r,k)}(x|\lambda)+T_{n}^{(r,k)}(x|\lambda).\\
\end{split}
\end{equation}

Applying $t$ on the both sides of Theorem\ref{thm3}, we get

\begin{equation}\label{40}
\begin{split}
&(n+1)T_{n}^{(r,k)}(x|\lambda)\\
&=nxT_{n-1}^{(r,k)}(x|\lambda)+T_{n}^{(r,k)}(x|\lambda)-rnT_{n-1}^{(r,k)}(x|\lambda)-\frac{rn\lambda}{1-\lambda}T_{n-1}^{(r+1,k)}(x|\lambda)\\
& -\frac{1}{n+1}\sum_{l=0}^{n+1}\binom{n+1}{l}B_{l}\{(n+1-l) T_{n-l}^{(r,k)}(x|\lambda)-(n+1-l)T_{n-l}^{(r,k-1)}(x|\lambda)\}\\
\end{split}
\end{equation}

Thus, by (\ref{40}), we have

\begin{equation}\label{41}
\begin{split}
&(n+1)T_{n}^{(r,k)}(x|\lambda)+n(r-\frac{1}{2}-x)T_{n-1}^{(r,k)}(x|\lambda)+\sum_{l=0}^{n-2}\binom{n}{l}B_{n-l} T_{l}^{(r,k)}(x|\lambda)\\
&=-\frac{r\lambda n}{1-\lambda}T_{n-1}^{(r+1,k)}(x|\lambda)+\sum_{l=0}^{n}\binom{n}{l}B_{n-l} T_{l}^{(r,k-1)}(x|\lambda).\\
\end{split}
\end{equation}

Therefore, by (\ref{41}), we obtain the following theorem.

\begin{theorem}\label{thm4}
For $r,k \in \mathbb{Z}, ~n \in \mathbb{Z}~ with~n\geq2$, we have
\begin{equation*}
\begin{split}
&(n+1)T_{n}^{(r,k)}(x|\lambda)+n(r-\frac{1}{2}-x)T_{n-1}^{(r,k)}(x|\lambda)+\sum_{l=0}^{n-2}\binom{n}{l}B_{n-l} T_{l}^{(r,k)}(x|\lambda)\\
&=-\frac{r\lambda n}{1-\lambda}T_{n-1}^{(r+1,k)}(x|\lambda)+\sum_{l=0}^{n}\binom{n}{l}B_{n-l} T_{l}^{(r,k-1)}(x|\lambda).\\
\end{split}
\end{equation*}
\end{theorem}

From  (\ref{14}) and  (\ref{26}), we note that

\begin{equation}\label{42}
\begin{split}
T_{n}^{(r,k)}(y|\lambda)&= \langle \left( \frac{1-\lambda}{e^{t}-\lambda} \right)^{r} \frac{Li_{k}(1-e^{-t})}{1-e^{-t}}e^{yt} \vert x^{n} \rangle\\
&= \langle \left( \frac{1-\lambda}{e^{t}-\lambda} \right)^{r} \frac{Li_{k}(1-e^{-t})}{1-e^{-t}}e^{yt} \vert x x^{n-1} \rangle.\\
\end{split}
\end{equation}

By (\ref{15}) and  (\ref{42}), we get

\begin{equation}\label{43}
\begin{split}
T_{n}^{(r,k)}(y|\lambda)&= \langle \partial_{t} \left(\left( \frac{1-\lambda}{e^{t}-\lambda} \right)^{r} \frac{Li_{k}(1-e^{-t})}{1-e^{-t}}e^{yt}\right) \vert x^{n-1} \rangle\\
&= \langle \left( \partial_{t}\left( \frac{1-\lambda}{e^{t}-\lambda} \right)^{r}\right) \frac{Li_{k}(1-e^{-t})}{1-e^{-t}}e^{yt} \vert x^{n-1} \rangle\\
&+ \langle \left( \frac{1-\lambda}{e^{t}-\lambda} \right)^{r} \left( \partial_{t}\frac{Li_{k}(1-e^{-t})}{1-e^{-t}}\right)e^{yt}  \vert x^{n-1} \rangle\\
&+ \langle \left( \frac{1-\lambda}{e^{t}-\lambda} \right)^{r}\frac{Li_{k}(1-e^{-t})}{1-e^{-t}} \partial_{t}e^{yt}  \vert x^{n-1} \rangle.\\
\end{split}
\end{equation}

Therefore, by (\ref{43}), we obtain the following theorem.

\begin{theorem}\label{thm5}
For $r,k \in \mathbb{Z}, ~n \geq1$, we have
\begin{equation*}
\begin{split}
&T_{n}^{(r,k)}(x|\lambda)=(x-r)T_{n-1}^{(r,k)}(x|\lambda)-\frac{r\lambda }{1-\lambda}T_{n-1}^{(r+1,k)}(x|\lambda)\\
&+\sum_{l=0}^{n-1}\left\{(-1)^{n-1-l}\binom{n-1}{l}\sum_{m=0}^{n-1-l}(-1)^{m}\frac{(m+1)!}{(m+2)^{k}}S_{2}(n-1-l,m)\right\}H_{l}^{(r)}(x-1|\lambda).\\
\end{split}
\end{equation*}
\end{theorem}

Now, we compute $\langle \left( \frac{1-\lambda}{e^{t}-\lambda} \right)^{r} Li_{k}(1-e^{-t}) \vert x^{n+1} \rangle$ in two different ways.\\

On the one hand,

\begin{equation}\label{44}
\begin{split}
& \langle \left( \frac{1-\lambda}{e^{t}-\lambda} \right)^{r} Li_{k}(1-e^{-t}) \vert x^{n+1} \rangle =\langle \left( \frac{1-\lambda}{e^{t}-\lambda} \right)^{r} \frac{Li_{k}(1-e^{-t})}{1-e^{-t}} \vert (1-e^{-t})x^{n+1} \rangle\\
&=\langle \left( \frac{1-\lambda}{e^{t}-\lambda} \right)^{r} \frac{Li_{k}(1-e^{-t})}{1-e^{-t}} \vert x^{n+1}-(x-1)^{n+1} \rangle\\
&=\sum_{m=0}^{n}\binom{n+1}{m}(-1)^{n-m}\langle \left( \frac{1-\lambda}{e^{t}-\lambda} \right)^{r} \frac{Li_{k}(1-e^{-t})}{1-e^{-t}} \vert x^{m} \rangle\\
&=\sum_{m=0}^{n}\binom{n+1}{m}(-1)^{n-m}\langle 1 \vert T_{m}^{(r,k)}(x|\lambda)\rangle\\
&=\sum_{m=0}^{n}\binom{n+1}{m}(-1)^{n-m} T_{m}^{(r,k)}(\lambda).\\
\end{split}
\end{equation}

On the other hand, we get

\begin{equation}\label{45}
\begin{split}
& \langle \left( \frac{1-\lambda}{e^{t}-\lambda} \right)^{r} Li_{k}(1-e^{-t}) \vert x^{n+1} \rangle =\langle Li_{k}(1-e^{-t}) \vert \left( \frac{1-\lambda}{e^{t}-\lambda} \right)^{r} x^{n+1} \rangle\\
&=\langle \int_{0}^{t}(Li_{k}(1-e^{-s}))'ds \vert H_{n+1}^{(r)}(x|\lambda)\rangle\\
&=\langle \int_{0}^{t}e^{-s}\frac{Li_{k}(1-e^{-s})}{(1-e^{-s})} ds \vert H_{n+1}^{(r)}(x|\lambda)\rangle\\
&=\sum_{l=0}^{n} \left( \sum_{m=0}^{l}\binom{l}{m}(-1)^{l-m}B_{m}^{(k-1)} \right) \frac{1}{l!}\langle \int_{0}^{t}s^{l}ds \vert H_{n+1}^{(r)}(x|\lambda)\rangle\\
&=\sum_{l=0}^{n} \sum_{m=0}^{l}\binom{l}{m}(-1)^{l-m}\frac{B_{m}^{(k-1)}}{(l+1)!} \langle t^{l+1} \vert H_{n+1}^{(r)}(x|\lambda)\rangle\\
&=\sum_{l=0}^{n} \sum_{m=0}^{l}\binom{l}{m}\binom{n+1}{l+1}(-1)^{l-m}B_{m}^{(k-1)}  H_{n-l}^{(r)}(\lambda).\\
\end{split}
\end{equation}

Therefore, by (\ref{44}) and (\ref{45}), we obtain the following theorem.

\begin{theorem}\label{thm6}
For $r,k \in \mathbb{Z}, ~n \in \mathbb{Z}_{\geq 0}$, we have
\begin{equation*}
\sum_{m=0}^{n}\binom{n+1}{m}(-1)^{n-m}T_{m}^{(r,k)}(\lambda)=\sum_{l=0}^{n} \sum_{m=0}^{l}(-1)^{l-m}\binom{l}{m}\binom{n+1}{l+1}B_{m}^{(k-1)}  H_{n-l}^{(r)}(\lambda).
\end{equation*}
\end{theorem}

Now, we consider the following two Sheffer sequences:

\begin{equation}\label{46}
\begin{split}
&T_{n}^{(r,k)}(x|\lambda) \sim \left( \left( \frac{e^{t}-\lambda}{1-\lambda} \right)^{r} \frac{1-e^{-t}}{Li_{k}(1-e^{-t})},t \right),\\
& ~~ \mathbb{B}^{(s)} \sim \left( \left( \frac{e^{t}-1}{t} \right)^{s} ,t \right),\\
\end{split}
\end{equation}
where $s \in \mathbb{Z}_{\geq 0}$, $r,k \in \mathbb{Z}$ and $\lambda \in \mathbb{C}$ with $\lambda \neq 1$.
Let us assume that

\begin{equation}\label{47}
T_{n}^{(r,k)}(x|\lambda) = \sum_{m=0}^{n}C_{n.m}\mathbb{B}^{(s)}_{m}(x).
\end{equation}

By (\ref{21}) and (\ref{47}), we get

\begin{equation}\label{48}
\begin{split}
C_{n,m}&=\frac{1}{m!} \langle \left( \frac{e^{t}-1}{t} \right)^{s} \left( \frac{1-\lambda}{e^{t}-\lambda} \right)^{r} \frac{Li_{k}(1-e^{-t})}{1-e^{-t}} t^{m} \vert x^{n} \rangle\\
&=\frac{1}{m!} \langle \left( \frac{e^{t}-1}{t} \right)^{s} \left( \frac{1-\lambda}{e^{t}-\lambda} \right)^{r} \frac{Li_{k}(1-e^{-t})}{1-e^{-t}} \vert t^{m}x^{n} \rangle\\
&=\binom{n}{m} \langle \left( \frac{e^{t}-1}{t} \right)^{s} \left( \frac{1-\lambda}{e^{t}-\lambda} \right)^{r} \frac{Li_{k}(1-e^{-t})}{1-e^{-t}} \vert x^{n-m} \rangle\\
&=\binom{n}{m} \sum_{l=0}^{n-m}\frac{s!}{(l+s)!}S_{2}(l+s,s)\langle  \left( \frac{1-\lambda}{e^{t}-\lambda} \right)^{r} \frac{Li_{k}(1-e^{-t})}{1-e^{-t}}  \vert t^{l}x^{n-m} \rangle\\
&=\binom{n}{m} \sum_{l=0}^{n-m}\frac{s!l!}{(l+s)!}\frac{(n-m)_{l}}{l!}S_{2}(l+s,s)\langle 1 \vert T_{n-m-l}^{(r,k)}(x|\lambda) \rangle \\
&=\binom{n}{m} \sum_{l=0}^{n-m}\frac{\binom{n-m}{l}}{\binom{s+l}{l}}S_{2}(l+s,s)T_{n-m-l}^{(r,k)}(\lambda).\\
\end{split}
\end{equation}

Therefore, by (\ref{47}) and (\ref{48}), we obtain the following thoerem.

\begin{theorem}\label{thm7}
For $r,k \in \mathbb{Z}, ~s \in \mathbb{Z}_{\geq 0}$, we have
\begin{equation*}
T_{n}^{(r,k)}(x|\lambda)=\sum_{m=0}^{n}\left\{ \binom{n}{m} \sum_{l=0}^{n-m}\frac{\binom{n-m}{l}}{\binom{s+l}{l}}S_{2}(l+s,s)T_{n-m-l}^{(r,k)}(\lambda) \right\} \mathbb{B}^{(s)}_{m}(x).
\end{equation*}
\end{theorem}

From (\ref{3}) and (\ref{22}), we note that

\begin{equation}\label{49}
\begin{split}
&T_{n}^{(r,k)}(x|\lambda) \sim \left( \left( \frac{e^{t}-\lambda}{1-\lambda} \right)^{r} \frac{1-e^{-t}}{Li_{k}(1-e^{-t}) },t \right),\\
& E_{n}^{(r,s)}(x) \sim \left( \left( \frac{e^{t}+1}{2} \right)^{s} ,t \right),\\
\end{split}
\end{equation}
where $r,k \in \mathbb{Z}, ~s \in \mathbb{Z}_{\geq 0}$.

By the same method, we get

\begin{equation}\label{50}
T_{n}^{(r,k)}(x|\lambda)=\frac{1}{2^{s}}\sum_{m=0}^{n}\left\{\binom{n}{m}\sum_{j=0}^{s}\binom{s}{j}T_{n-m}^{(r,k)}(j)\right\}E_{m}^{(s)}(x).
\end{equation}

From (\ref{1}) and (\ref{22}), we note that

\begin{equation}\label{51}
\begin{split}
&T_{n}^{(r,k)}(x|\lambda) \sim \left( \left( \frac{e^{t}-\lambda}{1-\lambda} \right)^{r} \frac{1-e^{-t}}{Li_{k}(1-e^{-t})},t \right),\\
& H_{n}^{(s)}(x| \mu ) \sim \left(\left( \frac{e^{t}-\mu}{1-\mu} \right)^{s} ,t \right),\\
\end{split}
\end{equation}
where $r,k \in \mathbb{Z}$,and $\lambda,~\mu \in \mathbb{C}$ with $\lambda \neq 1,~\mu \neq 1$, $ ~s \in \mathbb{Z}_{\geq 0}$.

Let us assume that
\begin{equation}\label{52}
T_{n}^{(r,k)}(x|\lambda)=\sum_{m=0}^{n}C_{n,m}H_{m}^{(s)}(x|\mu).
\end{equation}

By (\ref{21}) and (\ref{52}), we get

\begin{equation}\label{53}
\begin{split}
C_{n,m}&=\frac{1}{m!} \langle \left( \frac{e^{t}-\mu}{1-\mu} \right)^{s} \left( \frac{1-\lambda}{e^{t}-\lambda} \right)^{r} \frac{Li_{k}(1-e^{-t})}{1-e^{-t}} t^{m} \vert x^{n} \rangle\\
&=\frac{\binom{n}{m}}{(1-\mu)^{s}}\langle \left( e^{t}-\mu \right)^{s} \vert \left( \frac{1-\lambda}{e^{t}-\lambda} \right)^{r} \frac{Li_{k}(1-e^{-t})}{1-e^{-t}} x^{n-m} \rangle\\
&=\frac{\binom{n}{m}}{(1-\mu)^{s}}\sum_{j=0}^{s}\binom{s}{j}(-\mu)^{s-j}\langle e^{jt} \vert T_{n-m}^{(r,k)}(x|\lambda)\rangle\\
&=\frac{\binom{n}{m}}{(1-\mu)^{s}}\sum_{j=0}^{s}\binom{s}{j}(-\mu)^{s-j}T_{n-m}^{(r,k)}(j|\lambda).\\
\end{split}
\end{equation}

Therefore, by (\ref{52}) and (\ref{53}), we obtain the following theorem.

\begin{theorem}\label{thm8}
For $r,k \in \mathbb{Z}, ~s \in \mathbb{Z}_{\geq 0}$, we have
\begin{equation*}
T_{n}^{(r,k)}(x|\lambda)=\frac{1}{(1-\mu)^{s}}\sum_{m=0}^{n}\left\{\binom{n}{m}\sum_{j=0}^{s}\binom{s}{j}(-\mu)^{s-j}T_{n-m}^{(r,k)}(j|\lambda) \right\} H_{m}^{(s)}(x|\mu).
\end{equation*}
\end{theorem}

It is known that

\begin{equation}\label{54}
\begin{split}
&T_{n}^{(r,k)}(x|\lambda) \sim \left( \left( \frac{e^{t}-\lambda}{1-\lambda} \right)^{r} \frac{1-e^{-t}}{Li_{k}(1-e^{-t})},t \right),\\
& (x)_{n} \sim (1,~e^{t}-1).\\
\end{split}
\end{equation}

Let

\begin{equation}\label{55}
T_{n}^{(r,k)}(x|\lambda)=\sum_{m=0}^{n}C_{n,m}(x)_{m}.
\end{equation}

Then, by (\ref{21}) and (\ref{55}), we get

\begin{equation}\label{56}
\begin{split}
C_{n,m}&=\frac{1}{m!} \langle \left( \frac{1-\lambda}{e^{t}-\lambda} \right)^{r} \frac{Li_{k}(1-e^{-t})}{1-e^{-t}}(e^{t}-1)^{m} \vert x^{n} \rangle\\
&=\sum_{l=0}^{\infty} \frac{S_{2}(l+m,m)}{(l+m)!}\langle \left( \frac{1-\lambda}{e^{t}-\lambda} \right)^{r} \frac{Li_{k}(1-e^{-t})}{1-e^{-t}} \vert t^{m+l}x^{n} \rangle\\
&=\sum_{l=0}^{n-m} \frac{S_{2}(l+m,m)}{(l+m)!}(n)_{m+l}\langle 1 \vert \left( \frac{1-\lambda}{e^{t}-\lambda} \right)^{r} \frac{Li_{k}(1-e^{-t})}{1-e^{-t}}x^{n-m-l} \rangle\\
&=\sum_{l=0}^{n-m} \binom{n}{l+m} S_{2}(l+m,m)T_{n-m-l}^{(r,k)}(\lambda).\\
\end{split}
\end{equation}

Therefore, by (\ref{55}) and (\ref{56}), we  obtain the following theorem.

\begin{theorem}\label{thm9}
For $r,k \in \mathbb{Z}$, we have
\begin{equation*}
T_{n}^{(r,k)}(x|\lambda)=\sum_{m=0}^{n}\left\{ \sum_{l=0}^{n-m} \binom{n}{l+m} S_{2}(l+m,m)T_{n-m-l}^{(r,k)}(\lambda) \right\}(x)_{m}.
\end{equation*}
\end{theorem}

Finally, we consider the following two Sheffer sequences:

\begin{equation}\label{57}
\begin{split}
&T_{n}^{(r,k)}(x|\lambda) \sim \left( \left( \frac{e^{t}-\lambda}{1-\lambda} \right)^{r} \frac{1-e^{-t}}{Li_{k}(1-e^{-t})},t \right),\\
& x^{[n]} \sim (1,~1-e^{-t}),\\
\end{split}
\end{equation}
where $x^{[n]} = x(x+1) \cdots(x+n-1)$.

Let us assume that

\begin{equation}\label{58}
T_{n}^{(r,k)}(x|\lambda)=\sum_{m=0}^{n}C_{n,m}x^{[m]}.
\end{equation}

Then, by (\ref{21}) and (\ref{58}), we  get

\begin{equation}\label{59}
\begin{split}
C_{n,m}&=\frac{1}{m!} \langle \left( \frac{1-\lambda}{e^{t}-\lambda} \right)^{r} \frac{Li_{k}(1-e^{-t})}{1-e^{-t}}(1-e^{-t})^{m} \vert x^{n} \rangle\\
&=\sum_{l=0}^{\infty} \frac{(-1)^{l}S_{2}(l+m,m)}{(l+m)!}\langle \left( \frac{1-\lambda}{e^{t}-\lambda} \right)^{r} \frac{Li_{k}(1-e^{-t})}{1-e^{-t}} \vert t^{m+l}x^{n} \rangle\\
&=\sum_{l=0}^{n-m} \frac{(-1)^{l}S_{2}(l+m,m)}{(l+m)!}(n)_{m+l}\langle 1 \vert \left( \frac{1-\lambda}{e^{t}-\lambda} \right)^{r} \frac{Li_{k}(1-e^{-t})}{1-e^{-t}}x^{n-m-l} \rangle\\
&=\sum_{l=0}^{n-m} (-1)^{l}\binom{n}{l+m} S_{2}(l+m,m)T_{n-m-l}^{(r,k)}(\lambda).\\
\end{split}
\end{equation}

Therefore, by by (\ref{58}) and (\ref{59}), we  obtain the following theorem.

\begin{theorem}\label{thm10}
For $r,k \in \mathbb{Z},~n \geq 0$, we have
\begin{equation*}
T_{n}^{(r,k)}(x|\lambda)=\sum_{m=0}^{n}\{ \sum_{l=0}^{n-m}(-1)^{l} \binom{n}{l+m} S_{2}(l+m,m)T_{n-m-l}^{(r,k)}(\lambda) \}x^{[m]}.
\end{equation*}
\end{theorem}

%%%%%%%%%%%%%%%%%%%%%%%%%%%%%%%%%%%%%%%%%%%%%%%%%%%%%%%%%%%%%%%

\bigskip
ACKNOWLEDGEMENTS. This work was supported by the National Research Foundation of Korea(NRF) grant funded by the Korea government(MOE)\\
(No.2012R1A1A2003786 ).
\bigskip
%%%%%%%%%%%%%%%%%%%%%%%%%%%%%%%%%%%%%%%%%%%%%%%%%%%%%%%%%%%%%%

\end{document}